\documentclass[12pt]{amsart}

\usepackage{amsmath}
\usepackage{amsthm}
\usepackage[colorlinks]{hyperref}
\usepackage{enumerate}

\newcommand{\cb}{\mathcal B}

\newcommand{\gtg}{G \times \widehat{G}}

\allowdisplaybreaks

\theoremstyle{plain}
\newtheorem{thm}{Theorem}[section]
\newtheorem{defn}[thm]{Definition}
\newtheorem{lem}[thm]{Lemma}
\newtheorem{prop}[thm]{Proposition}

\newtheorem{rem}[thm]{Remark}
\newtheorem{exam}[thm]{Example}

\title{Vector-valued properties of the Weyl transform}

\author[R. Singhal]{Ritika Singhal$^*$} 
\author[N. S. Kumar]{N. Shravan Kumar}
\thanks{ \newline
	 Ritika Singhal$^{1,*}$, N. Shravan Kumar$^{2}$\\
	Department of Mathematics, Indian Institute of Technology, Delhi,\\
	Hauz Khas, New Delhi-110016,  India\\
	Email$^1$: ritikasinghal1120@gmail.com \\
 Email$^2$: shravankumar.nageswaran@gmail.com\\
	 * Corresponding author}
	
\begin{document}
	
	\begin{abstract}
		In this paper, we introduce and study the concept of Weyl type and cotype w.r.t a locally compact abelian group. We also study the Weyl transform of vector measures and functions which are integrable w.r.t a vector measure. Later, we also introduce and study the convolution of functions from $L^p$-spaces associated to a vector measure.
	\end{abstract}
	
	
	
	\maketitle

	\section{Introduction}
	
  The notion of Weyl transform was introduced by Hermann Weyl in his book on {\it the theory of groups and quantum mechanics} \cite{Weyl}. A large number of researchers from both Mathematics and Physics have studied this transform, including its applications to partial differential equations and quantum physics.

  The study of vector-valued functions has attracted increasing attention in recent decades. Most of the classical problems of the theory of functions can be studied in a vector-valued setting. It allows to combine independent results of the theory of functions based on a common approach. In \cite{garcia}, Garc\'{i}a et al. showed the importance of vector-valued analysis via maximum Hardy-Littlewood operator and theorems of Littlewood-Paley type (\cite{stein}, Chapter 4).  The authors of \cite{garcia} left the study of vector-valued functions open, clearing that further investigations in this field are an important and interesting problem and may lead to new, unexpected results. Recently, properties of pseudo-differential operators with symbols from operator-valued spaces have been considered in \cite{XX}. Following the path, we studied the Weyl transform of vector-valued functions. Properties like the Hausdorff-Young inequality and the Plancheral theorem introduced for the Fourier transform also held for the Weyl transform. We investigated the properties further in a vector-valued setting and showed that the analogue of the above properties still holds for the Weyl transform of the vector-valued functions, but the Reimann-Lebesgue lemma fails for the Weyl transform of vector measures (Example \ref{RLMF1}). At the same time, we also provide a sufficient condition for a vector measure to satisfy the Riemann-Lebesgue lemma (Proposition \ref{SCRLL}).

	The notion of Fourier type w.r.t a locally compact abelian group was introduced by J. Peetre \cite{Pe} for $\mathbb{R}$ and by Milman \cite{M} for general locally compact abelian groups. This notion was extended to compact groups by Garc\'{i}a-Cuerva and Parcet \cite{GP} and by Hun Hee Lee \cite{HHL} for general locally compact unimodular groups.  Following the ideas of the articles mentioned above, we introduce and study the Weyl type and cotype of operators in the framework of operator spaces. 

	On the other hand, in the past few years, the well-known results of the classical abelian harmonic analysis have been extended to the setup of vector measures. In \cite{CFNP}, the author studied the Fourier transform and a vector-valued convolution for (weakly) $\nu$-integrable functions. Blasco in \cite{B} analyzed the Fourier transform of vector measures $\nu$ as well as the convolution between scalar and vector-valued regular measures for a compact abelian group. Very recently, M. Kumar and N. Kumar in \cite{MkNsk} studied the analogue of the classical Young's inequality and other properties on a compact group.
	
	The main aim of Section 4 of this paper is to study the Weyl transform and related properties w.r.t a vector measure. Here, we define and study the Weyl transform of functions which are integrable (weakly) w.r.t a vector measure. We also study the Weyl transform of vector-valued measures in the same section. We want to mention here that this is the first instance where a study of vector measures, from the point of view of harmonic analysis, has been carried out on non-compact groups. Finally, Section 5 defines the concept of twisted convolution of functions w.r.t vector measures. Two different kinds of convolution are defined and are also shown to be equivalent. 
	We shall begin with preliminaries on the Weyl transform and operator spaces required in the sequel.
	
	\section{Preliminaries}
	\subsection{Weyl transform and twisted convolution} 
	Let $G$ be a locally compact abelian group with $\widehat{G}$ as its dual group.  We shall denote by $m_{G\times\widehat{G}},$ the Haar measure on the product group $G\times\widehat{G}.$  For $1\leq p \leq \infty$, we shall denote by $L^p(G\times\widehat{G})$, the classical $L^p$-spaces and by $C_c(\gtg)$ the space of compactly supported functions on $G\times\widehat{G}$. Throughout, $p'$ will denote the conjugate index of $p$.
	 	
 	The \textit{Weyl transform}, denoted $W,$ is defined as a $\mathcal{B}\left(L^{2}(G)\right)$-valued integral on $C_{c}(G \times \widehat{G})$, given by 
 	$$(W(f) \varphi)(y)=\int_{G \times \widehat{G}}f(x, \chi) \rho((x, \chi))(\varphi)(y) dx d\chi, f \in C_c(G \times \widehat{G})$$ 
 	 where $\rho(x,\chi)$ is motivated from the \textit{Schr\"odinger representation} of the abstract Heisenberg group $\mathbb{H}(G)$ on $L^2(G)$ and defined as $$\rho((x, \chi))(\varphi)(y)= \chi(y) \varphi(x y), \varphi \in L^{2}(G).$$
 	
	For $f, g \in C_{c}(G \times \widehat{G}),$ the {\it twisted convolution}, denoted $f \times g,$ is defined as $$f \times g(x, \chi):=\int_{G \times \widehat{G}} f\left(x x'^{-1}, \chi \overline{\chi'}\right) g\left(x^{\prime}, \chi^{\prime}\right)\chi'(xx'^{-1}) d m_{G\times\widehat{G}}( x',  \chi').$$
	 Here are some integrability properties of the twisted convolution.
	 See \cite{RS}.
	\begin{prop}\label{TCProp1}
		Suppose $1 \leq p \leq \infty, f \in L^{1}(G \times \widehat{G})$ and $g \in L^{p}(G \times \widehat{G})$.
		\begin{enumerate}[i)]
			\item Then $f \times g \in L^{p}(G \times \widehat{G})$ and $\|f \times g\|_{p} \leq\|f\|_{1}\|g\|_{p}$.
			\item Also $g \times f \in L^{p}(G \times \widehat{G})$ and $\|g \times f\|_{p} \leq\|f\|_{1}\|g\|_{p}$.
			\item When $p=\infty,$ both $f \times g$ and $g \times f$ are continuous.
		\end{enumerate}
	\end{prop}
	\begin{prop}\label{TCp4}
		Suppose $f \in L^{p}(G \times \widehat{G})$ and $g \in L^{q}(G \times \widehat{G})$. 
		\begin{enumerate}[i)]
			\item For $q=p'$, $f \times g \in C_{0}(G \times \widehat{G})$ and $\|f \times g\|_{\infty} \leq\|f\|_{p}\|g\|_{q}$.
			\item If $p=q=2,$ then $f \times g \in L^{2}(G \times \widehat{G})$ and $\|f \times g\|_{2} \leq\|f\|_{2}\|g\|_{2}$.
			\item If $\frac{1}{p} + \frac{1}{q} = \frac{1}{r}+1,$ then $f \times g \in L^{r}(G \times \widehat{G})$ and $\|f \times g\|_{r} \leq\|f\|_{p}\|g\|_{q}$. 
		\end{enumerate}
	\end{prop}
	Let $\mathcal{H}$ be a Hilbert space and then for  $1\leq p<\infty$,  $\mathcal{B}_p(\mathcal{H})$ denotes  Schatten $p$-class of $\mathcal{H}$ and for $p=\infty$ it represents the space of all compact operators on $\mathcal{H}.$
	As a consequence of the definition of the twisted convolution and the above integrability properties, we have the following. 
	
	\begin{thm} \label{WTProp}\cite{RS}
		\begin{enumerate}[i)]
			
			\item For $p\in\{1,2\}$, the space $L^p(G \times \widehat{ G})$ is a Banach algebra with pointwise addition and twisted convolution as addition and multiplication, respectively.
			\item (Plancherel Theorem) The Weyl transform is an isometric Banach algebra isomorphism between $L^2(G\times\widehat{G})$ and $\mathcal{B}_2(L^2(G)).$
			\item (Riemann-Lebesque Lemma) The Weyl transform is a bounded operator and a Banach algebra homomorphism from $L^1(G \times \widehat{G})$ to $\mathcal{B}_\infty\left(L^{2}(G)\right).$
			\item (Hausdorff-Young inequality) If $1<p<2,$ then the Weyl transform maps $L^p(G\times\widehat{G})$ into $\mathcal{B}_{p'}(L^2(G)).$
		\end{enumerate}
	\end{thm}
	
$W_{G,p}$ will denote the Weyl transform on $L^p(\gtg)$. We will skip the subscript when the domain is obvious.	For more on Weyl transform and twisted convolution, see \cite{RS}.
	
	\subsection{Operator spaces} \label{os}
	An \textit{operator space} is a closed subspace of $\mathcal{B}(\mathcal{H})$ for a Hilbert space $\mathcal{H}.$ Let $X$ and $Y$ be operator spaces. Let $M_n(X)$ denote the space of all $n\times n$ matrices with entries from $X.$ 	
	A natural norm on $M_n(X)$ can be defined using the natural identification of  $M_n (\mathcal{B}(\mathcal H))$ with $\mathcal{B} (\oplus_{i}^n \mathcal{H})$. Let $T:X\rightarrow Y$ be a linear map, then for $n\in\mathbb{N},$ the $n^{th}$-amplification of $T,$ denoted $T^{(n)},$ is a linear map from $M_n(X)$ to $M_n(Y)$ given by $T^{(n)}([x_{ij}])=[T(x_{ij})].$ The map $T$ is \textit{completely bounded} if $$\|T\|_{cb}:=\sup\{\|T^{(n)}\|:n\in\mathbb{N}\}<\infty.$$
	We shall denote by $\mathcal{CB}(X)$, the space of all completely bounded maps from $X$ to $\mathbb{C}$, equipped with the norm $\|\cdot\|_{cb}$. Also, $T$ is a \textit{complete isometry}(or \textit{contraction}) if for each $n\in\mathbb{N},$ $T^{(n)}$ is an isometry(or contraction).
	Let $X\otimes_{\min}Y$ and $X \widehat{\otimes} Y$ denote the \textit{minimal tensor product} and operator space \textit{projective tensor product} of $X$ and $Y$, respectively. It is worth mentioning here that inclusion $X\widehat{\otimes}Y\hookrightarrow X\otimes_{\min} Y$ holds completely contractively.
	
	Let $\{X, Y\}$ be a compatible couple of Banach spaces in the sense of complex interpolation. If $X_\theta,\ 0<\theta<1,$ denotes the interpolation space $(X,Y)_\theta,$ then, by \cite[Pg. 53]{P1}, $X_\theta$ is also an operator space. For $1\leq p \leq \infty,$ the space $\mathcal{B}_p[\mathcal{H};X]$ is defined as 
	$$\mathcal{B}_p[\mathcal{H};X]:=(\mathcal{B}_\infty\otimes_{\min}X,\mathcal{B}_1\widehat{\otimes}X)_{1/p}.$$ 
	We shall denote by $S^n_p(X)$ the space $\mathcal{B}_p[\ell^2(n);X]$ where $\ell^2(n)=(\mathbb{C}^n,\|\cdot\|_2)$. See \cite{ER, P1} for more on operator spaces. The following theorem summarizes the properties of these interpolated spaces. 
	
	\begin{thm}[\cite{P2}, Chapter 1]\label{pisier}\mbox{ } Let $1\leq p,p_0,p_1\leq \infty.$ Then we have the following:
		\begin{enumerate}[i)]
			\item The $cb$-norm of  $T:X\rightarrow Y$ is equal to $\underset{n\geq 1}{\sup}\|I_{M_n}\otimes T\|_{\mathcal{B}(S^n_p(X),S^n_p(Y))}.$
			\item (Duality) The dual of $\mathcal{B}_p[\mathcal{H};X]$ is completely isometrically isomorphic to $\mathcal{B}_{p^\prime}[\mathcal{H};X^\ast].$ 
			\item (Fubini type theorem) For $n\geq 1,$ the following are complete isometry:
			\begin{align*}
				S_p^n\left(L^p(G\times\widehat{G};X)\right)\cong & L^p\left(G\times\widehat{G},S_p^n(X)\right) \\
				S^n_p\left(\mathcal{B}_p[\mathcal{H};X]\right)\cong & \mathcal{B}_p[\mathcal{H};S^n_p(X)].
			\end{align*} 
			\item (Interpolation) If $0<\theta<1$, $1\leq p_0,p_1\leq \infty$ and $\frac{1}{p_\theta}=\frac{1-\theta}{p_0}+\frac{\theta}{p_1},$ then we have $$\mathcal{B}_{p_\theta}[\mathcal{H};X]=\left(\mathcal{B}_{p_0}[\mathcal{H};X],\mathcal{B}_{p_1}[\mathcal{H};X]\right)_\theta$$ completely isometrically.
		\end{enumerate}
	\end{thm}

	\subsection{Vector measures} \label{vm}
		Let $A$ denote a complex Banach space and $(\Omega,\mathcal{A})$ denote a measurable space. The space of all $\sigma$-additive, $A$-valued vector measures on $\Omega$ will be denoted by $M(\Omega,A)$. Let $\nu $ be in $M(\Omega,A)$. Let $A^*$ denote the dual of $A$ and $B_{A}$ be the closed unit ball in $A.$ 
  For $B\in\mathcal{A}$, the \textit{semivariation of $\nu$ }  is given by $\|\nu\|(B)=\underset{x^*\in B_{A^\ast}}\sup |\langle\nu,x^\ast\rangle|(B)$.
		
	A function $f: \Omega \to \mathbb{C}$ is said to be $\nu$-weakly integrable if
	$$\|f\|_\nu=\underset{x^\ast\in B_{A^\ast}}{\sup}\int_\Omega |f|\,d|\langle\nu,x^\ast\rangle|< \infty.$$ Let $L^1_w(\Omega,\nu)$ denote the set of all $\nu$-weakly integrable functions. If for each  $B\in\mathcal{A},$ there exists an unique element,  denoted by $\int_Bfd\nu$, in $A$ such that $\int_B f\,d\langle\nu,x^\ast\rangle=\langle \int_Bfd\nu,x^\ast\rangle,\ x^\ast\in A^*$, then $f$ is said to be \textit{$\nu$-integrable}.
	The space of all $\nu$-integrable functions will be denoted by $L^1(\Omega,\nu)$. For $1 < p < \infty ,L^p_w(\Omega, \nu)$ and $L^p(\Omega, \nu)$ consist of the functions  $f: \Omega \to \mathbb{C}$ such that $f^p \in L^1_w(\Omega, \nu)$ and $f^p \in L^1(\Omega, \nu)$, respectively. Both the spaces are Banach if equipped with the norm $\|f\|_{\nu,p}=\|f^p\|^{1/p}_\nu$.  The space $L^\infty(\Omega, \nu)=L^\infty_w(\Omega, \nu)$ will denote the space of all $\nu$-a.e. bounded functions, which is also a Banach space with the essential supremum norm $\|\cdot\|_{\nu,\infty}$. 
	
	 Consider the group $G \times \widehat{G}$ with Borel $\sigma$-algebra, $\mathfrak{B}(G\times \widehat{G})$. The space $\mathcal{M}(G \times\widehat{G},A)$ consists of all  regular $\sigma$-finite vector measures on $A$. A vector measure $\nu$ is a.c. w.r.t a non-negative scalar measure $\mu$ if $\underset{\mu(B)\rightarrow 0}{\lim}\nu(B)=0,\ B\in\mathcal{A}.$ The space
	$M_{ac}(G\times\widehat{G},A) \subset \mathcal{M}(G\times\widehat{G},A) $ consists of all $\nu \in M(G\times\widehat{G},A)$ such that $\nu$ is a.c. w.r.t the Haar measure on $G \times \widehat{G}$.

		If $f \in L^1(G \times \widehat{G},\nu)$, then $\nu_f(U):=\int_U f d \nu $ defines a vector measure on $\mathfrak{B}(G \times \widehat{G})$ with $\|\nu_f\|=\|f\|_\nu$. For $f \in L^1_w(G \times \widehat{G},\nu)$, $\nu_f$ is a $X^{**} $-valued measure on $\mathfrak{B}(G \times \widehat{G})$ given by $\langle \nu_f(U),x^* \rangle = \int_U f d \langle \nu, x^* \rangle , U \in \mathfrak{B}(G \times \widehat{G}), x^* \in B_{X^*}.$
	The \textit{$p$-semivariation} of $\nu$ w.r.t $m_{G \times \widehat{G}}$,
	for $1 \leq p<\infty$, is given by
	$$\|\nu\|_{p, m_{G\times \widehat{G}}}=\sup \left\{\bigg \|\sum_{B \in \rho} \alpha_{B} \nu(B) \bigg \|_{A}: \begin{array}{c}
		\rho \text { is the finite partition of $A$ with } \\ \sum_{B \in \rho} \alpha_{B} \chi_{B} \in B_{L^{p^\prime}(G\times \widehat{G})}
	\end{array} \right\}$$
	
	and for $p=\infty,\|\nu\|_{\infty, m_{G\times \widehat{G}}}=\underset{m_{G\times \widehat{G}}(B)>0}{\sup} \frac{\|\nu(B)\|}{m_{G\times\widehat{G}}(B)}$ . $ M_{p}(G\times \widehat{G};A)$ will denote the space of all $A$-valued $\sigma$-additive vector measures with finite $p$-semivariation. For $1 \leq p \leq \infty$, $P_{p}(G\times \widehat{G};A)$ denotes the closure of the space of all $A$-valued simple functions on $G\times \widehat{G}$ in $M_{p}(G\times \widehat{G};A)$, equipped with the norm
	$$
	\|\phi\|_{P_{p}(G\times \widehat{G};A)}=\|\nu_\phi\|_{p,m_{G \times \widehat{G}}}=\sup _{x^{\ast} \in B_{A^\ast}}\left\|\left\langle\phi, x^{\ast}\right\rangle\right\|_{p}.
	$$
	The space $C_c(G\times \widehat{G};A)$, consisting of all $A$ -valued continuous functions on $G\times\widehat{G}$ having compact support, is dense in $P_{p}(G\times \widehat{G};A) $, for $ 1 \leq p<\infty$ and closed in $P_{\infty}(G;A).$ The space $L^p(G \times \widehat{G};A)$ stands for the space of all $A$-valued $p$-Bochner integrable functions. 
	 We shall refer the interested readers to \cite{D, DU, HNVW} for more details on vector measures.
	
	\begin{lem} \label{Contain}
		If $\nu \in M( \Omega; A) ,$ then for $1<p \leq \infty$, we have $L^p_w(\Omega, \nu) \subset L^1_w(\Omega, \nu).$
	\end{lem}
	\begin{proof}
		Without any loss of generality, we assume $\|\nu\|=1.$  If $ f \in  L^p(\Omega, \nu)$, 
		then $f^p \in L^1_w(\nu)$, i.e $\underset{x^\ast\in B_{A^\ast}}{\sup}\int_\Omega |f|^pd|\langle\nu,x^\ast\rangle| < \infty.$
		Hence, using Holder's inequality,
		\begin{eqnarray*}
			\int_\Omega |f|\,d|\langle\nu,x^\ast\rangle| & \leq \left( \int_\Omega |f|^p\,d|\langle\nu,x^\ast\rangle| \right)^{1/p} \left(
			|\langle \nu,x^* \rangle| (\Omega)\right)^{1/q}.
		\end{eqnarray*}
		Taking supremum over $x^* \in A^*$ in both sides, we get $\|f\|_\nu \leq \|f\|_{\nu,p} $. \qedhere
		
	\end{proof}
	The following result is of independent interest, and a version of this is proved in \cite{B}.
	\begin{lem}\label{dense}
		Let $\nu \in \mathcal{M}(G \times \widehat{G},A)$. Then $C_c(G\times\widehat{G})$  is dense in $L^p(G\times\widehat{G},\nu),$ for $1\leq p<\infty.$
	\end{lem}
	\begin{proof}
		As simple functions are dense in $L^p(G\times\widehat{G}, \nu),$ it suffices to show that for any Borel set $A$ of $G\times\widehat{G},$ the function $\chi_A$ can be approximated by a function from $C_c(G\times\widehat{G})$ in the $L^p(G\times\widehat{G}, \nu)$-norm. Let $A\in\mathfrak{B}(G\times\widehat{G})$ and $\epsilon>0.$ By assumption, $\nu$ is regular and therefore, by definition, there exists a compact set $K$ and an open set $U$ such that $K\subset A\subset U$ and $\|\nu\|(U\setminus K)<\epsilon^p.$ Now, by Urysohn's lemma, $\exists$ $f\in C_c(G\times\widehat{G})$ such that $0\leq f\leq 1,$ $f$ is identically $1$ on $K$ and $\mbox{supp}(f)\subset U.$ 
		
		We now claim that this $f$ approximates $\chi_A.$ In fact,
		\begin{eqnarray*}
			\|\chi_A-f\|_{L^p(G\times\widehat{G},\nu)} &=& \underset{x^\ast\in B_{A^*}}{\sup}\|\chi_A-f\|_{L^p(G\times\widehat{G}, |\langle \nu,x^\ast \rangle|)} \\ &=& \underset{x^\ast\in B_{A^*}}{\sup}\left( \int_{U\setminus K} |\chi_A-f|^p\ d|\langle \nu,x^\ast \rangle| \right)^{1/p} \\ &\leq& (\|\nu\|(U\setminus K))^{1/p}<\epsilon.
		\end{eqnarray*}
		Hence, the proof.
	\end{proof}
	
	\subsection{Tensor integrability} \label{ti}
	In this subsection, we provide the basics related to tensor integrability. The concepts presented here are already contained in \cite{S}. 
	Let $X$ and $Y$ be operator spaces and  $\nu \in M(G \times \widehat{G};X)$.
 
 For a $\nu$-measurable function $f: \Omega \rightarrow Y$, we shall denote by $N(f)$, the quantity
	$$\sup_{x^ * \in B_{X^ *}} \int \|f\|\ d|\langle \nu, x^ * \rangle |.$$
	A $\nu $-measurable $Y$-valued function $f$ is said to be {\it $\otimes_{\min }$-integrable} if there exists a sequence $\left(\phi_{n}\right)$ of  $Y$-valued simple functions such that $\lim _{n} N\left(\phi_{n}-f\right)=0$ and $\int_{A} f d \nu$ represents the limit of the
sequence $\left(\int_{A} \phi_{n} d \nu \right)$ in  $Y \otimes_{\min } X$.  Let $L^{1}(\Omega, \nu, Y, X)$ denote the space of all $\otimes_{\text {min}}$-integrable functions on $\Omega$ which becomes a Banach space when equipped with $N(\cdot)$ norm.
	\begin{thm}
		Let $f: \Omega \rightarrow Y$ be a $\nu$-measurable function. Then $f \in L^{1}(\Omega, \nu, Y, X)$ if and only if $\|f\| \in L^{1}(\Omega, \nu)$.
	\end{thm}
	
	\begin{prop}\label{OprandInt_Commu}
		If $f \in L^{1}(\Omega,\nu, Y, X),$ then for $y^{*} \in Y^{*}$ and $T \in \mathcal{C B}(X)$,
		$$
		\left(y^{*} \otimes T\right)\left(\int_{A} f d \nu \right)=\int_{A}\left\langle f, y^{*}\right\rangle d(T \circ \nu), A \in \Sigma .$$
	\end{prop}
	A weakly $\|\nu\|$-measurable function $f: \Omega \rightarrow Y$ is said to have a {\it generalized  weak $\otimes$-integral} w.r.t $\nu$ if $\langle f, y^\ast \rangle \in L^1_w(\Omega,\nu),\forall \  y^* \in Y^*.$ We shall denote by gen-$L_{w}^{1}(\Omega,\nu, Y, X)$ the space consisting of $Y$-valued generalized weak $\otimes$-integrable functions on $\Omega$. It is a normed linear space when equipped with the norm 
	$$ N_{w}(f)=\sup _{y^{*} \in B_{Y^{*}}}\left\|\left\langle f, y^{*}\right\rangle\right\|_{\nu}. $$
	
	With all these basics on tensor integrability, we now show that the spaces $L^1(\Omega,\nu),$ $L^1_w(\Omega,\nu)$ and $M(\Omega;x)$ are all operator spaces if the Banach space $X$ is also an operator space. So, in the next three results, $X$ is a complete operator space.
	\begin{thm}\label{OSSL1}
		The following spaces are isomorphic via the mapping $
		\left[f_{i j}\right] \mapsto \widetilde{f}$,  where $ \widetilde{f}(\cdot)=\left[f_{i j}(\cdot)\right].		$
		\begin{enumerate}[i)]
			\item $M_{n}\left(L^{1}(\Omega,\nu)\right) $and $ L^{1}\left(\Omega,\nu, M_{n}, X\right)$;
			\item $M_{n}\left(L_{w}^{1}(\Omega,\nu)\right)  $ and gen-$L_{w}^{1}\left(\Omega,\nu, M_{n}, X\right)$;
			\item $M_{n}(M(\Omega;X)) $ and $ M\left(\Omega, M_{n}(X)\right)$.
		\end{enumerate}
		The matrix norm arising from the above identification gives an operator space structure on $L^1(\Omega,\nu),$ $L^1_w(\Omega,\nu)$ and $M(\Omega; X)$.
		
	\end{thm}
	\section{Weyl transform of vector-valued functions}
	In this section, we study the vector-valued analogue of the Hausdorff-Young inequality for the Weyl transform associated to a locally compact abelian group. In order to establish this result, we introduce the concepts of Weyl type and Weyl cotype. Throughout this section, let $X$  be an operator space and  $1\leq p\leq 2$.
	
	\subsection{Weyl type w.r.t a locally compact abelian group}
	We shall begin by introducing the concept of the Weyl type.
	\begin{defn}
		We say that a linear map $T:X\rightarrow Y$ is of $G$-Weyl type $p$ if $W_{G,p}\otimes T:L^p(G\times\widehat{G})\otimes X\rightarrow\mathcal{B}_{p^\prime}(L^2(G))\otimes Y$ extends to a completely bounded map, still denoted $W_{G,p}\otimes T,$ from $L^p(G\times\widehat{G};X)$ into $\mathcal{B}_{p^\prime}[L^2(G);Y].$
	\end{defn}
	
	Our first result shows that if an operator has Weyl type for some $p$, it is completely bounded.
	\begin{lem}\label{WTpCB}
		Let $T: X\rightarrow Y$ be a linear map. If $T$ has $G$-Weyl type $p,$ then $T$ is completely bounded and $$\|T\|_{cb}\leq \|W_{G,p}\|^{-1}\|W_{G,p}\otimes T\|_{cb}.$$
	\end{lem}
	\begin{proof}
		Suppose that $T$ has $G$-Weyl type $p,$ i.e., the map $W_{G,p}\otimes T:L^p(G\times\widehat{G};X)\rightarrow\mathcal{B}_{p^\prime}[L^2(G);Y]$ is completely bounded. This in turn, by $i)$ of Theorem \ref{pisier}, is equivalent to saying that $$I_{M_n}\otimes(W_{G,p}\otimes T):S^n_{p^\prime}(L^p(G\times\widehat{G};X))\rightarrow S^n_{p^\prime}(\mathcal{B}_{p^\prime}[L^2(G);Y])$$ is uniformly bounded (in $n$). Note that, by $iii)$ of Theorem \ref{pisier}, $S^n_{p^\prime}(\mathcal{B}_{p^\prime}[L^2(G);Y])$ is isometrically isomorphic to $\mathcal{B}_{p^\prime}[L^2(G);S^n_{p^\prime}(Y)].$ Further, if $f\in L^p(G\times\widehat{G})$ and $[x_{ij}]\in M_n(X),$ then note that $f\otimes[x_{ij}]$ is mapped to $W_{G,p}(f)\otimes[T(x_{ij})].$ Thus, \begin{eqnarray*}
			\|W_{G,p}(f)\|_{\cb_{p^{\prime}}(L^2(G))}\|[T(x_{ij})]\|_{S^n_{p^\prime}(Y)} &=& \|W_{G,p}(f)\otimes [T(x_{ij})]\|_{\cb_{p^{\prime}}[L^2(G);S^n_{p^\prime}(Y)]} \\ &=& \|(W_{G,p}\otimes T^{(n)}) (f\otimes [x_{ij}])\|_{\cb_{p^{\prime}}[L^2(G);S^n_{p^\prime}(Y)]} \\ &\leq& \|W_{G,p}\otimes T^{(n)}\|\|f\otimes [x_{ij}]\|_{S^n_{p^\prime}(L^p(G\times\widehat{G};X))} \\ &\leq& \|W_{G,p}\otimes T\|_{cb}\|f\otimes [x_{ij}]\|_{L^p(G\times\widehat{G},S^n_{p^\prime}(X))} \\ &=& \|W_{G,p}\otimes T\|_{cb}\|f\|_{L^p(G\times\widehat{G})} \|[x_{ij}]\|_{S^n_{p^\prime}(X)}.
		\end{eqnarray*}
		Since $f\in L^p(G\times\widehat{G})$ is arbitrary, we are done.
	\end{proof}
	The following lemma shows that the converse of the above lemma is true for the case $p=1.$ This is the vector-valued version of the Riemann-Lebesgue lemma for the Weyl transform.
	\begin{lem}\label{WT1}
		If $T:X\rightarrow Y$ is a completely bounded linear map, then $T$ has $G$-Weyl type $1.$ Moreover, $$\|W_{G,1}\otimes T\|_{cb}=\|T\|_{cb}.$$
	\end{lem}
	\begin{proof}
		Suppose that $T$ is completely bounded. Without loss of generality, we can assume that $T$ is completely contractive. By Theorem \ref{WTProp}, the Weyl transform map $W_{G,1}:L^1(G\times\widehat{G})\rightarrow\mathcal{B}_\infty(L^2(G))$ is bounded. In fact, it is contractive. We shall equip $L^1(G\times\widehat{G})$ with the natural maximal operator space structure so that, by \cite[Pg. 49]{ER}, $W_{G,1}$ is completely bounded. Thus, by \cite[Corollary 7.1.3]{ER}, the map $W_{G,1}\otimes T:L^1(G\times\widehat{G})\otimes X\rightarrow \mathcal{B}_\infty(L^2(G))\otimes Y$ extends to a complete contraction from $L^1(G\times \widehat{G};X)\rightarrow\mathcal{B}_\infty(L^2(G))\widehat{\otimes}Y.$ On the other hand, by \cite[Pg. 142]{ER}, the space $\mathcal{B}_\infty(L^2(G))\widehat{\otimes}Y$ sits inside $\mathcal{B}_\infty[L^2(G);Y]$ completely contractively. 
	\end{proof}
	We say that an operator space $X$ has {\it $G$-Weyl type $p$} if the identity operator on $X$ has $G$-Weyl type $p.$ The following few lemmas give some examples for operator spaces with $G$-Weyl type $p$.
	\begin{lem}\label{WTLp}
		Let $(\Omega,\mathcal{A},\mu)$ be a measure space. Then the space $L^p(\Omega)$ has $G$-Weyl type $p.$ Similarly, if $\mathcal{H}$ is any Hilbert space, then $\mathcal{B}_p(\mathcal{H})$ has $G$-Weyl type $p.$
	\end{lem}
	\begin{proof}
		For the case of $p=2,$ this is a consequence of $ii)$ of Theorem \ref{WTProp} and \cite[Theorem 1.1]{P3}. The general case then follows from complex interpolation.
	\end{proof}
	Again, by using complex interpolation, we can obtain the following lemma.
	\begin{lem}\label{WTpq}
		Let $1\leq p<q\leq 2.$ If $X$ has $G$-Weyl type $q$, then $X$ has $G$-Weyl type $p.$
	\end{lem}
	The following lemma is about the relation between the Weyl type w.r.t. product of abelian groups and the Weyl type w.r.t. the individual groups.
	\begin{lem}\label{WTProduct}
		Let $G$ and $H$ be locally compact abelian groups. An operator $T:X\rightarrow Y$ has $(G\times H)$-Weyl type $p$ if and only if $W_{H,p}\otimes T$ has $G$-Weyl type $p.$
	\end{lem}
	\begin{proof}
		The proof of this is a consequence of the fact that $W_{G\times H,p}$ coincides with $W_{G,p}\otimes W_{H,p}.$
	\end{proof}
	Again, we have the following lemma by the same reasoning given above and by Lemma \ref{WTpCB}.
	\begin{lem}\label{WTProduct1}
		Let $G$ and $H$ be locally compact abelian groups, $T:X\rightarrow Y$ and $S:Y\rightarrow Z$ be completely bounded linear operators between operator spaces such that $T$ has $G$-Weyl type $p$ and $S$ has $H$-Weyl type $p.$ Then
		\begin{enumerate}[i)]
		  \item $\|W_{G,p}\otimes(S\circ T)\|_{cb}\leq \|T\|_{cb}\|W_{G,p}\otimes S\|_{cb}.$
			\item $\|W_{H,p}\otimes T\|_{cb}\leq\|W_{G,p}\|^{-1}\|W_{G\times H,p}\otimes T\|_{cb}.$
			\item $\|W_{G\times H,p}\otimes(S\circ T)\|_{cb}\leq\|W_{G,p}\otimes T\|_{cb}\|W_{H,p}\otimes S\|_{cb}.$
		\end{enumerate}
	\end{lem}
	
	\subsection{Weyl cotype w.r.t. a locally compact abelian group}
	We now introduce the concept of the Weyl cotype.
	\begin{defn}
		We say that a linear map $T:X\rightarrow Y$ is of $G$-Weyl cotype $p^\prime$ if $W_{G,p}^{-1}\otimes T:\mathcal{B}_p(L^2(G))\otimes X\rightarrow L^{p^\prime}(G\times\widehat{G})\otimes Y$ extends to a completely bounded map from $\mathcal{B}_p[L^2(G);X]$ into $L^{p^\prime}(G\times\widehat{G};Y).$
	\end{defn}
	The following is an analogue of Lemma \ref{WTpCB} for the case of the Weyl cotype.
	\begin{lem}\label{WcTinfty}
		Suppose that $T:X\rightarrow Y$ is a linear map and has $G$-Weyl cotype $p^\prime.$ Then $T$ is completely bounded and $$\|T\|_{cb}\leq \|W_{G,p}^{-1}\|^{-1}\|W_{G,p}^{-1}\otimes T\|_{cb}.$$
	\end{lem}
	\begin{proof}
		Since the proof of this is the same as Lemma \ref{WTpCB}, we omit the proof.
	\end{proof}
	The following lemma is the converse of the above result when $p^\prime=\infty.$
	\begin{lem}
		If $T:X\rightarrow Y$ is a completely bounded linear map, then $T$ has $G$-Weyl cotype $\infty.$ Moreover, $$\|W_{G,p}^{-1}\otimes T\|_{cb}=\|T\|_{cb}.$$
	\end{lem}
	\begin{proof}
		As the proof of this is the same as the proof of Lemma \ref{WT1}, we omit the proof.
	\end{proof}
	We say that an operator space $X$ has {\it $G$-Weyl cotype $p^\prime$} if the identity operator on $X$ has $G$-Weyl cotype $p^\prime.$ The next two lemmas give some examples for operator spaces with $G$-Weyl cotype $p^\prime$ and are analogues of Lemma \ref{WTLp}, and Lemma \ref{WTpq}, respectively. As the proofs are similar, we omit them.
	\begin{lem}
		Let $(\Omega,\mathcal{A},\mu)$ be a measure space. Then the space $L^p(\Omega)$ has $G$-Weyl cotype $p^\prime.$ Similarly, if $\mathcal{H}$ is any Hilbert space, then $\mathcal{B}_p(\mathcal{H})$ has $G$-Weyl cotype $p^\prime.$
	\end{lem}
	\begin{lem}
		Let $1\leq p<q\leq 2.$ If $X$ has $G$-Weyl cotype $q^\prime$, then $X$ has $G$-Weyl cotype $p^\prime.$
	\end{lem}
	Similarly, we can give the Weyl cotype result for Lemma \ref{WTProduct} and Lemma \ref{WTProduct1}.
	Here is the main result of this section. This is the vector-valued version of the duality of the Hausdorff-Young inequality for the Weyl transform.
	\begin{thm}
		Let $T: X\rightarrow Y$ be a linear map and let $T^*$ denote its Banach space adjoint.
		\begin{enumerate}[i)]
			\item The map $T$ has $G$-Weyl type $p$ if and only if $T^*$ has $G$-Weyl cotype $p^\prime.$
			\item The map $T$ has $G$-Weyl cotype $p^\prime$ if and only if $T^*$ has $G$-Weyl type $p.$
		\end{enumerate}
	\end{thm}
	\begin{proof}
		Since the proof of $ii)$ is similar to the proof of $i)$, we omit its proof. Again, it is plain that in order to prove $i)$, it is enough to prove only one side as the proof of the other side follows similar lines. Thus, we shall now prove the only if part of $i)$.
		
		Let $[S_{ij}]\in S^n_{p^\prime}(\mathcal{B}_p[L^2(G);Y^\ast]).$ Let $\epsilon>0.$ By $ii)$ of Theorem \ref{pisier}, there exists $[f_{ij}]\in S^n_p(L^p(G\times\widehat{G};X))$ such that $\|[f_{ij}]\|_{S^n_p(L^p(G\times\widehat{G};X))}=1$ and $$\|(W_{G,p}^{-1}\otimes T^\ast)([S_{ij}])\|_{S^n_{p^\prime}(L^{p^\prime}(G\times\widehat{G};X^\ast))}\leq (1+\epsilon)|\langle [f_{ij}], (W_{G,p}^{-1}\otimes T^\ast)([S_{ij}]) \rangle|.$$ Note that $W_{G,p}^{-1}\otimes T^\ast=(W_{G,p}\otimes T)^\ast.$ Thus, 
		\begin{eqnarray*}
			& & \|\left(W_{G,p}^{-1}\otimes T^\ast\right)([S_{ij}])\|_{S^n_{p^\prime}(L^{p^\prime}(G\times\widehat{G};X^\ast))} \leq  (1+\epsilon) \left| \langle [f_{ij}], \left(W_{G,p}^{-1}\otimes T^\ast\right)([S_{ij}])  \rangle \right| \\ &=& (1+\epsilon) \left| \langle [f_{ij}], \left(W_{G,p} \otimes T\right)^\ast([S_{ij}])  \rangle \right| = (1+\epsilon) \left| \langle \left(W_{G,p} \otimes T\right)\left([f_{ij}]\right), ([S_{ij}])  \rangle \right| \\ &\leq& (1+\epsilon) \|W_{G,p}\otimes T\|_{cb} \|[f_{ij}]\|_{S^n_p(L^p(G\times\widehat{G};X))}\|[S_{ij}]\|_{S^n_{p^\prime}(\mathcal{B}_p[L^2(G);Y^\ast])} \\ &=& (1+\epsilon) \|W_{G,p}\otimes T\|_{cb} \|[S_{ij}]\|_{S^n_{p^\prime}(\mathcal{B}_p[L^2(G);Y^\ast])}.
		\end{eqnarray*}
		As $\epsilon>0$ is arbitrary, letting $\epsilon\rightarrow 0,$ we get the desired inequality.
	\end{proof}

	\section{Weyl transform on $L^1(G\times\widehat{G},\nu),$ $L^1_w(G\times\widehat{G},\nu)$ and $M(G\times\widehat{G};X)$} 
	In this section, we study the Weyl transform of $\nu$-integrable functions. We also study the Weyl transform of vector measures. We assume $\nu \in M(G \times \widehat{G}; X)$ for an operator space $X$. We show that the Weyl transform is a completely bounded operator. We also provide an example to show the failure of an analogue of the Riemann-Lebesgue lemma. For each $x^* \in X^* , h_{x^*} $ represents the Radon-Nikodym derivative $\frac{d\langle \nu ,x^*\rangle }{dm_{G \times \widehat{G}}}$.
	\subsection{Weyl transform on $L^1(G\times\widehat{G},\nu)$}\label{transform}
	We first define the notion of the Weyl transform of functions in $L^1(G\times\widehat{G},\nu).$ Observe that, if $f\in L^1(G\times\widehat{G},\nu)$, then $f\rho \in L^1(G\times\widehat{G}, \nu,\mathcal{B}(L^2(G));X).$ In fact, this helps us define the following:
	
	\begin{defn} \label{def1}
		The Weyl transform of a function $f\in L^1(G\times\widehat{G},\nu)$ w.r.t. the vector measure $\nu$ is defined by $$W^\nu(f)=\int_{G\times\widehat{G}}f(x,\chi)\rho(x,\chi)\ d\nu(x,\chi).$$
	\end{defn}
	Note that, for $f\in L^1(G\times\widehat{G},\nu),$ $W^\nu(f)\in \mathcal{B}(L^2(G))\otimes_{min}X = \mathcal{B}[L^2(G);X].$
	
	\begin{exam}\label{WTEx}
		Let $U$ be an open subset of $G\times\widehat{G}$ having compact closure. Let $1\leq p<\infty$ and let $T:L^p(\overline{U})\rightarrow X$ be a completely bounded operator. Consider the vector measure $\nu:\mathfrak{B}(G\times\widehat{G})\rightarrow X$ given by $\nu(A)=T(\chi_{A\cap U}).$ Note that $L^p(\overline{U})\subset L^1(G\times\widehat{G},\nu).$ Using the density of simple functions, it can be shown that for $A\in\mathfrak{B}(G\times\widehat{G})$ and $f\in L^p(G\times\widehat{G},\nu),$ $$\int_A f\ d\nu=T(f\chi_{A\cap U}).$$ Let $g=\sum_{i=1}^n y_i \otimes \chi_{A_i}$ be any $Y$-valued simple function. Then
		\begin{eqnarray*} \int_A g d\nu&=& \sum_{i=1}^n y_i \otimes \nu(A \cap A_i)  \\&=&\sum_{i=1}^n (Id_Y \otimes T )\left(y_i \otimes\chi_{A \cap A_i \cap U}\right) \\&=& (Id_Y \otimes T) \left(g \cdot \chi_{A \cap U}\right). \end{eqnarray*}
		Hence for $f\in L^1(G\times\widehat{G},\nu)$ and $g\in L^1(\mathfrak{B}(G\times\widehat{G}),\nu,Y;X) $,
		$$\int_A(fg)d \nu=(Id_Y \otimes T)((f \chi_{A\cap U}) g).$$
		Thus, if $f\in L^1(G\times\widehat{G}, \nu  ),$ then $$W^\nu(f)=(Id_{\mathcal{B}(L^2(G))}\otimes T)((f\chi_{U})\rho).$$
	\end{exam}
	\begin{lem}
		The mapping $f\mapsto W^\nu(f)$ from $L^1(G\times\widehat{G},\nu)$ to $\mathcal{B}(L^2(G))\otimes_{\min}X$ is completely bounded.
	\end{lem}
	\begin{proof}
		The proof follows from  \cite[Corollary 3]{CSb}.
	\end{proof}
	We now define the Weyl transform of $\nu$-weakly integrable functions. For $f\in L^1_w(G \times \widehat{G},\nu),$ $$\int_{G\times\widehat{G}}|f(x,\chi)|\ d|\langle \nu,x^\ast \rangle|(x,\chi)=\int_{G\times\widehat{G}} |f(x,\chi)h_{x^\ast}(x,\chi)|\ dm_{G\times\widehat{G}}(x,\chi),$$ i.e., $fh_{x^\ast}\in L^1(G\times\widehat{G})$ $\forall\ f\in L^1_w(G\times\widehat{G},\nu).$ This observation leads to the following definition.
	\begin{defn}\label{def2}
		Let $\nu\in M_{ac}(G\times\widehat{G};X).$ If $f\in L^p_w(G\times\widehat{G}, \nu),$ then by Lemma \ref{Contain},  $f\in L^1_w(G\times\widehat{G}, \nu).$ Also, for each $x^*\in X^*,$ $fh_{x^*}\in L^1(G\times\widehat{G})$. Thus, the Weyl transform of $f\in L^p_w(G\times\widehat{G},\nu)$ is defined by $$W^\nu(f)(x^*):=W(fh_{x^*}),\ x^*\in X^*.$$
	\end{defn}
 \begin{rem}
    It is well known that $X\otimes Y$ can be seen as a subset of $\mathcal{B}(Y^*, X)$. Under this identification, for $f \in L^1(G \times \widehat{G}, \nu)$, Definition \ref{def1} and Definition \ref{def2} coincide.
\end{rem}

	\begin{lem}
		Let $\nu\in M_{ac}(G\times\widehat{G};X).$ 
		\begin{enumerate}[i)]
			\item For $f\in L^1_w(G\times\widehat{G},\nu),$ the mapping $x^\ast\mapsto W^\nu(f)(x^*)$ from $X^\ast$ to $\mathcal{B}(L^2(G))$ is bounded.
			\item The mapping $f\mapsto W^\nu(f)$ from $L^1_w(G\times\widehat{G},\nu)$ to $\mathcal{B}(X^\ast,\mathcal{B}(L^2(G)))$ is completely bounded.
		\end{enumerate}
	\end{lem}
	\begin{proof}
		Note that $i)$ is a routine check while $ii)$ is a consequence of Theorem \ref{OSSL1}.
	\end{proof}
	
	We now prove the uniqueness theorem for the Weyl transform.
	\begin{thm}[Uniqueness theorem]
		Let $\nu\in M_{ac}(G\times\widehat{G};X)$ and $f\in L^1_w(G\times\widehat{G},\nu).$ If $W^\nu(f)=0,$ then $f=0$ $\nu$-a.e..
	\end{thm}
	\begin{proof}
		Let $f\in L^1_w(G\times\widehat{G},\nu)$ such that $W^\nu(f)=0.$ Then, by definition, $W^\nu(f)(x^*)=0$ for every $x^*\in X^*,$ i.e., $W(fh_{x^*})=0$ for every $x^*\in X^*.$ This implies by the uniqueness theorem for the Weyl transform that $fh_{x^*}=0\ m_{G\times\widehat{G}}$-a.e., i.e., $\exists$ a Borel set $A\subset G\times\widehat{G}$ such that $fh_{x^*}=0$ on $(G\times\widehat{G})\setminus A$ and $m_{G\times\widehat{G}}(A)=0.$ As shown in \cite[Theorem 4.9]{MkNsk}, one can show that $f=0\ \nu$-a.e..
	\end{proof} 
 
	A natural question that one would like to ask here is about an analogue of the Riemann-Lebesgue lemma.  To be more precise, if $f\in L^1(G\times\widehat{G},\nu)$, then will it imply that $W^\nu(f)\in\mathcal{B}_\infty[L^2(G);X].$
 The following example shows that the analogue of the Riemann-Lebesgue lemma fails.
	\begin{exam}\label{RLMF1}
		Let $G$ be an infinite group. Let $X=L^1(G\times\widehat{G})$ and let $U\subset G\times\widehat{G}$ be an open set having compact closure. Observe that $L^1(\overline{U})\subset L^1(G\times\widehat{G}).$ Let $T$ denote the  inclusion of $L^1(\overline{U})$ inside $L^1(G\times\widehat{G}).$ Let $\nu$ be the vector measure given in Example \ref{WTEx}.
  Let $0\neq f\in L^1(\overline{U}).$ Then \begin{align*}
		  W^\nu(f) = (f\chi_U)\rho \in & L^1(G\times\widehat{G},\mathcal{B}(L^2(G))) \\ \cong & \mathcal{B}(L^2(G))\widehat{\otimes} L^1(G\times\widehat{G})\subset\mathcal{B}[L^2(G);L^1(G\times\widehat{G})].
		\end{align*}
  Thus, if $W^\nu(f)\in \mathcal{B}_\infty[L^2(G);L^1(G\times\widehat{G})]$ for all $f\in L^1(G\times\widehat{ G}),$ then this would imply that $(f\chi_U)\rho\in L^1(G\times\widehat{ G},\mathcal{B}_\infty(L^2(G))).$ This in turn will imply that $f(x,\chi)\rho(x,\chi)$ is a compact operator for a.e. $(x,\chi)\in U,$ which is a contradiction. 
	\end{exam} 
	
	\subsection{Weyl transform on $M(G\times\widehat{G};X)$}
	
	We now define the Weyl transform of vector measures.
	\begin{defn}
		The Weyl transform of a measure $\nu\in M(G\times\widehat{G};X)$ is defined by $$W(\nu)=\int_{G\times\widehat{G}}\rho(x,\chi)\ d\nu(x,\chi).$$
	\end{defn}
	If $f\in L^1(G\times\widehat{G},\nu),$ then we shall denote by $\nu_f$ the corresponding vector measure, given by $$\nu_f(A)=\int_A f(x,\chi)\ d\nu(x,\chi).$$ The following lemma shows that the Weyl transform of $f$ and $\nu_f$ coincides.
	\begin{lem}\label{SPWTVM}
		\mbox{ }
		\begin{enumerate}[i)]
			\item If $f\in L^1(G\times\widehat{G},\nu),$ then $W^\nu(f)=W(\nu_f).$
			\item If $\nu\in M_{ac}(G\times\widehat{G};X)$ and $f\in L^1_w(G\times\widehat{G},\nu),$ then for $x^\ast\in X^*$ and $\varphi\in\mathcal{B}(L^2(G))^*,$ we have $$\langle W^\nu(f)(x^\ast),\varphi \rangle=\langle W^\nu(f),\varphi\otimes x^\ast \rangle$$
		\end{enumerate}
	\end{lem}
	\begin{proof}
		$i)$ is a consequence of the definition of $\nu_f$ while $ii)$ is a simple calculation.
	\end{proof}
	\begin{lem}
		The mapping $\nu\mapsto W(\nu)$ from $M(G\times\widehat{G};X)$ to $\mathcal{B}[L^2(G);X]$ is completely bounded.
	\end{lem}
	\begin{proof}
		The proof of this is a consequence of Theorem \ref{OSSL1}.
	\end{proof}
	We now prove the uniqueness theorem, assuming that the vector measure is a.c. w.r.t. $m_{G\times\widehat{G}}.$
	\begin{prop}[Uniqueness theorem]
		Let $\nu\in M_{ac}(G\times\widehat{G};X)$ be such that $W(\nu)=0.$ Then $\nu=0.$
	\end{prop}
	\begin{proof}
		Let $\varphi\in\mathcal{B}(L^2(G))^*$ and let $x^\ast\in X^*.$ Then, by $ii)$ of Lemma \ref{SPWTVM} \begin{eqnarray*}
			\langle W(\nu),\varphi\otimes x^\ast \rangle &=& \langle W(\nu)(x^\ast),\varphi \rangle \\ &=& \left\langle \int_{G\times\widehat{G}}\rho(x,\chi)\ d\langle \nu,x^\ast \rangle, \varphi \right\rangle \\ &=& \int_{G\times\widehat{G}} \langle \rho(x,\chi), \varphi \rangle\ d\langle \nu,x^\ast \rangle \\ &=& \int_{G\times\widehat{G}} \langle \rho(x,\chi), \varphi \rangle\ h_{x^\ast}(x,\chi)\ dm_{G\times\widehat{G}}(x,\chi) \\ &=& \left\langle \int_{G\times\widehat{G}} \rho(x,\chi)\ h_{x^\ast}(x,\chi)\ dm_{G\times\widehat{G}}(x,\chi) , \varphi \right\rangle = \langle W(h_{x^\ast}),\varphi \rangle.
		\end{eqnarray*}
		By hypothesis $W(\nu)=0.$ This implies that $\langle W(\nu),\varphi\otimes x^\ast \rangle = 0$ for all $\varphi\in\mathcal{B}(L^2(G))^*$ and for all $x^\ast\in X^*.$ By the above calculations, this is equivalent to saying that $\langle W(h_{x^\ast}),\varphi \rangle=0$ for all $\varphi\in\mathcal{B}(L^2(G))^*$ and for all $x^\ast\in X^*.$ Thus, $W(h_{x^\ast})=0$ for all $x^\ast\in X^*.$ By the uniqueness theorem for the Weyl transform, we have $h_{x^\ast}=0\, m_{G\times\widehat{G}}$-a.e. for all $x^\ast\in X^*.$ 
		
		Now, let $A$ be a Borel subset of $G\times\widehat{G}$ and let $x^\ast\in X^*.$ Then, $$|\langle \nu,x^\ast \rangle|(A)=\int_A\ d|\langle \nu,x^\ast \rangle| = \int_A |h_{x^\ast}|\ dm_{G\times\widehat{G}}=0.$$ Since $x^\ast\in X^*$ is arbitrary, it follows that $\nu=0.$
	\end{proof}
	Again, we provide an example to show the failure of the Riemann-Lebesgue lemma for the Weyl transform of vector measures.
	\begin{exam}
		Let $G$ be an infinite group. Let $X,$ $T$ and $\nu$ be as in Example \ref{RLMF1}. Note that $\nu\in M_{ac}(G\times\widehat{G},L^1(G\times\widehat{G})).$ Further, observe that $W(\nu)=\chi_U\rho$ and the same reasoning given in Example \ref{RLMF1} will tell us that $\nu$ does not satisfy the Riemann-Lebesgue property.
	\end{exam}
	Our next result gives a sufficient condition for a vector measure to satisfy the Riemann-Lebesgue Lemma.
	\begin{prop}\label{SCRLL}
		Let $\nu\in M_{ac}(G\times\widehat{G};X)$ be a measure of bounded variation. Also, let $X$ satisfy the Radon-Nikodym property w.r.t $(G\times\widehat{G}, \mathfrak{B}(G\times\widehat{G}),m_{G\times\widehat{G}})$ i.e. there exists $f\in L^1(G\times\widehat{G};X)$  such that $d\nu=f\,dm_{G \times \widehat{G}}.$ Then $W(\nu)\in\mathcal{B}_\infty[L^2(G);X].$ 
	\end{prop}
	\begin{proof}
		Since  $d\nu = fdm_{G\times\widehat{G}}$, we have $W(\nu)=W(f).$ Thus by Lemma \ref{WT1}, $W(\nu)\in\mathcal{B}_\infty[L^2(G);X].$
	\end{proof}
		\section{Twisted convolution} 
	In this section, we study the twisted convolution of $\nu$-weakly integrable functions. Here, we introduce and study two kinds of convolution and show that both definitions coincide. We also study the twisted convolution of a scalar measure and a vector measure. Let $X$ be an operator space and  $\nu \in M(G \times \widehat{G}; X)$.

\subsection{Twisted convolution associated  to a vector measure}\label{convo}
	We shall begin with the definition of twisted convolution of an $m_{G\times\widehat{G}}$-integrable function and a $ \nu$-weakly integrable function. Throughout this subsection, we assume $\nu \in M_{ac}(G \times \widehat{G};X)$ and $h_{x^*} $ represents the Radon-Nikodym derivative $\frac{d\langle \nu ,x^*\rangle }{dm_{G \times \widehat{G}}}$.
	\begin{defn}\label{TCvm1}
		Let $1\leq p\leq\infty.$ The twisted convolution of functions $f\in C_c(G\times\widehat{G})$ and $g\in L^1_w(G\times\widehat{G}, \nu)$ is defined as $$f\times_\nu g(x^*)=f\times (gh_{x^*}),\ x^*\in X^*.$$
	\end{defn}
	The following is an analogue of $i)$ of Proposition \ref{TCProp1}, and the proof is also a consequence of the same.
	\begin{lem}\label{p2}
		Let $1\leq p\leq\infty.$ If $f\in L^p(G\times\widehat{G})$ and $g\in L^1_w(G\times\widehat{G}, \nu),$ then $$f\times_\nu g\in\mathcal{B}(X^*,L^p(G\times\widehat{G}))$$ with $$\|f\times_\nu g\|_{\mathcal{B}(X^*,L^p(G\times\widehat{G}))}\leq \|f\|_p\|g\|_\nu.$$
	\end{lem}
	\begin{proof}
		Let $f\in L^p(G\times\widehat{G})$ and $g\in L^1_w(G\times\widehat{G}, \nu).$ Then, for $x^*\in B_{ X^*},$ we have
		\begin{align*}
			\|f\times_\nu g(x^*)\|_p =& \|f\times(gh_{x^*})\|_p\leq \|f\|_p\|gh_{x^*}\|_1 \leq \|f\|_p\|g\|_\nu.\qedhere
		\end{align*}
	\end{proof}

	The following is again a consequence of Proposition \ref{TCProp1} and \ref{Contain} whose proof follows the same lines as in Lemma \ref{p2}, and hence we omit it.
	\begin{lem}
		Let $1\leq p\leq\infty.$ If $f\in L^1(G\times\widehat{G})$ and $g\in L^p_w(G\times\widehat{G},\nu),$ then $f\times_\nu g\in \mathcal{B}(X^*,L^p(G\times\widehat{G})).$
	\end{lem}
	Our next lemma is a consequence of Proposition \ref{TCp4} and \ref{Contain}.
	\begin{lem}Let $1\leq p,q\leq\infty.$
		\begin{enumerate}[i)]
			\item  If  $f\in L^p(G\times\widehat{G})$ and $g\in L^{p'}_w(G\times\widehat{G},\nu),$ then $f \times_\nu g \in\mathcal{B}(X^*, C_0(G \times \widehat{G}))$ .
			\item If $f\in L^2(G\times\widehat{G})$ and $g\in L^2_w(G\times\widehat{G},\nu),$ then $f\times_\nu g\in \mathcal{B}(X^*,L^2(G\times\widehat{G}))$.
			\item If $\frac{1}{p}+\frac{1}{q}=\frac{1}{r}+1,$ $f\in L^p(G\times\widehat{G})$ and $g\in L^q_w(G\times\widehat{G},\nu),$ then $f\times_\nu g\in \mathcal{B}(X^*,L^r(G\times\widehat{G}))$ .
		\end{enumerate}
	\end{lem}
	Now, we have a lemma, which is a consequence of Lemma \ref{WTProp}.
	\begin{lem}
		For $p=1,2$, if $f\in L^p(G\times\widehat{G})$ and $g\in L^p_w(G\times\widehat{G}, \nu),$ then $$W(f\times_\nu g(x^*))=W(f)W^\nu(g)(x^*)\ \forall x^*\in X^*.$$
	\end{lem}
	\begin{proof}
		Let $f\in L^p(G\times\widehat{G})$ and $g\in L^p_w(G\times\widehat{G}, \nu).$ Then, for $x^*\in X^*,$ we have,
		\begin{align*}
			W(f\times_\nu g(x^*)) = & W(f\times g(h_{x^*})) \\ =& W(f)W(g(h_{x^*}))=W(f)W^\nu(g)(x^*).\qedhere
		\end{align*}
	\end{proof}

	\subsection{Vector-valued twisted convolution} \label{tcavm}
We now define vector-valued twisted convolution for any two measurable functions.
	\begin{defn}\label{TCvm2}
		We define the vector-valued twisted convolution, w.r.t $\nu,$ of two measurable functions $f$ and $g$, denoted $f \times ^{\nu} g$, as
		$$
		f \times ^{\nu} g(x,\chi)=\int_{G\times\widehat{G}}f\left(x x'^{-1}, \chi \overline{\chi'}\right) g\left(x^{\prime}, \chi^{\prime}\right)\chi'(xx'^{-1}) d \nu(x',\chi')
		$$
		provided that the mapping $(x',\chi ') \mapsto f\left(x x'^{-1}, \chi \overline{\chi'}\right) g\left(x^{\prime}, \chi^{\prime}\right)\chi'(xx'^{-1})$ belongs to $L^{1}(G \times \widehat{G},\nu)$ for $m_{G\times\widehat{G}}$ -a.e. $(x,\chi) \in G \times \widehat{G}$.
	\end{defn}
	The following lemma shows that Definition \ref{TCvm1} and Definition \ref{TCvm2} are equivalent.
	\begin{lem}\label{r1}
		Let $\nu \in M_{a c}(G \times \widehat{G};X) .$ If $f \in L^{p}(G \times \widehat{G}), 1 \leq p<\infty$ and $g \in L_{w}^{1}(G \times \widehat{G},\nu)$ are such that the mapping  $(x',\chi ') \mapsto f\left(x x'^{-1}, \chi \overline{\chi'}\right) g\left(x^{\prime}, \chi^{\prime}\right)\chi'(xx'^{-1})$ belongs to $ L^{1}(G \times \widehat{G},\nu)$ for $m_{G\times\widehat{G}}$ -a.e. $(x,\chi) \in G  \times \widehat{G} $, then for $x^{*} \in X^{*}$, we have $f \times _{\nu} g\left(x^{*}\right)=\left\langle f \times ^{\nu} g, x^{*}\right\rangle .$
	\end{lem}
	\begin{proof}
		Since $\nu \in  M_{a c}(G \times \widehat{G};X),$ for a given $x^*$ in $X^*,$ let $h_{x^*} \in L^1(G \times \widehat{G})$ be such that $d \langle \nu , x^* \rangle =h_{x^*}dm_{G\times\widehat{G}}.$ Then for $x^*$ in $X^*,$ we have
		\begin{align*}
			f \times _{\nu} g\left(x^{*}\right) =& f \times (gh_{x^*}) \\
			 =& \int_{G \times \widehat{G}} f\left(x x'^{-1}, \chi \overline{\chi'}\right) (gh_{x^*})\left(x^{\prime}, \chi^{\prime}\right)\chi'(xx'^{-1}) d m_{G\times\widehat{G}}( x',  \chi') \\
			  =& \int_{G \times \widehat{G}} f\left(x x'^{-1}, \chi \overline{\chi'}\right) g\left(x^{\prime}, \chi^{\prime}\right)\chi'(xx'^{-1})h_{x^*} ( x',  \chi')d m_{G\times\widehat{G}}(x',\chi') \\
			   =& \int_{G \times \widehat{G}} f\left(x x'^{-1}, \chi \overline{\chi'}\right) g\left(x^{\prime}, \chi^{\prime}\right)\chi'(xx'^{-1})d \langle \nu,x^* \rangle ( x',  \chi')\\
			   =&\left\langle f \times ^{\nu} g, x^{*}\right\rangle .  \qedhere	
		\end{align*}
	\end{proof}
	Before we proceed to the next result, here is the definition of Dunford's integrability.  A function $f:G\times\widehat{G}\rightarrow X$ is said to be \textit{Dunford integrable} if $\langle f,x^\ast \rangle\in L^1(G\times\widehat{G})\forall x^\ast\in X^\ast.$ We shall denote by $L^1_w(G\times\widehat{G};X)$ the space of Dunford integrable functions equipped with the norm $$\|f\|_{L^1_w(G\times\widehat{G};X)}=\underset{x^\ast\in B_{X^\ast}}{\sup}\|\langle f,x^\ast\rangle\|_{L^1(\gtg)}.$$
	
	A Dunford integrable function $f: G \times \widehat{G} \to X$ is said to be \textit{Pettis integrable }if for each $A \in \mathfrak{B}(G \times \widehat{G})$, there exists  $x_A \in X$ such that $\int_A \langle f, x^* \rangle dm_{G \times \widehat{G}}= \langle x_A,x^* \rangle$. The vector $x_A$ is unique and denoted by $(P)\int_A f dm_{G \times \widehat{G}}.$

	\begin{thm}
		Let $\nu \in M_{a c}(G\times \widehat{G};X) .$ If $f \in L^{1}(G\times \widehat{G})$ and $g \in L_{w}^{1}(G\times \widehat{G},\nu)$ are such that the mappings $(x',\chi ') \mapsto f(xx'^{-1},\chi\overline{\chi'})g(x',\chi ')\chi'(xx'^{-1}) $ are in $L^{1}(G\times \widehat{G},\nu)$ for $m_{G\times\widehat{G}}$ -a.e. $(x,\chi) \in G \times \widehat{G},$ then $f \times^{\nu} g$ is Dunford integrable with $\left\|f \times ^{\nu} g\right\|_{L_w^{1}(G\times \widehat{G};X) } \leq\|f\|_{1}\|g\|_{\nu} .$ In particular, if $g \in L^1(G \times \widehat{G},\nu)$, then $f \times^\nu g$ is Pettis integrable.	
	\end{thm}
	\begin{proof}
		Let $x^* \in B_{X^*}$ and let $h_{x^\ast} \in L^1(G \times \widehat{G})$ be such that $ d \langle \nu , x^\ast \rangle = h_{x^\ast} dm_{G\times\widehat{G}}.$ Define $\psi_{x^*}(x,\chi)=\langle f \times^\nu g(x,\chi), x^* \rangle,$ for $m_{G\times\widehat{G}}$-a.e. $(x,\chi) \in G \times \widehat{G}.$ Then $\psi_{x^*}(x,\chi)=f \times gh_{x^*}(x,\chi)$.
		Hence $\psi_{x^*}$ is measurable. Now, by Fubini's theorem, one can show that 
		\begin{align*}
			\int_{G\times \widehat{G}}\left|\left\langle f \times ^{\nu} g(x,\chi), x^{*}\right\rangle\right| d m_{G\times\widehat{G}}(x,\chi) &\leq\|f\|_{L^{1}(G\times \widehat{G})}\|g\|_{\nu},
		\end{align*}
		thus proving the Dunford integrability.
		
		Now, consider the mapping $(x',\chi') \mapsto \int_A  f\left(x x'^{-1}, \chi \overline{\chi'}\right) \chi'(x)  d m_{G \times \widehat{G}}(x,\chi)$. Since $|\chi'(x)|$ and $f \in L^1(G\times \widehat{G})$, the above mapping is bounded and measurable. Also, since $g \in L^1(G \times \widehat{G},\nu)$ and $|\chi'(x')|=1$, by \cite[Theorem 8]{dCFFMN}, the integral $$\int_{G \times \widehat{G}}g\left(x^{\prime}, \chi^{\prime}\right)\overline{\chi'(x')}d  \nu ( x',  \chi')$$ is well defined. 
		We claim $(P)\int_A f\times^\nu g dm_{G \times \widehat{G}}=x_A$ where $x_A$ is given by $$x_A= \int_{G \times \widehat{G}} \left( \int_A  f\left(x x'^{-1}, \chi \overline{\chi'}\right) \chi'(x) d m_{G \times \widehat{G}}(x,\chi) \right )g\left(x^{\prime}, \chi^{\prime}\right)\overline{\chi'(x')}d  \nu ( x',  \chi')\in X.$$ For $x^*\in X^*,$ we have, 
		\begin{align*} & \int_A \langle f \times ^ \nu g,x^*\rangle dm_{G \times \widehat{G}} \\ =& \int_A \left \langle  \int_{G \times \widehat{G}} f\left(x x'^{-1}, \chi \overline{\chi'}\right) g\left(x^{\prime}, \chi^{\prime}\right)\chi'(x) \overline{\chi'(x')} d \nu( x',  \chi'),x^* \right \rangle d m_{G \times \widehat{G}}(x,\chi) \\=& \int_A  \int_{G \times \widehat{G}} f\left(x x'^{-1}, \chi \overline{\chi'}\right) g\left(x^{\prime}, \chi^{\prime}\right)\chi'(x) \overline{\chi'(x')} d \langle \nu,x^* \rangle ( x',  \chi') d m_{G \times \widehat{G}}(x,\chi) \\
		=& \int_{G \times \widehat{G}} \left( \int_A  f\left(x x'^{-1}, \chi \overline{\chi'}\right) \chi'(x)  d m_{G \times \widehat{G}}(x,\chi) \right )g\left(x^{\prime}, \chi^{\prime}\right)\overline{\chi'(x')}d \langle \nu,x^* \rangle ( x',  \chi') \\
		=&  \langle x_A,x^* \rangle,\end{align*} thus proving our claim.
		\end{proof}
	
	\subsection{Twisted convolution of a vector measure} \label{TCVM}
	Finally, we define the twisted convolution of a scalar measure and a vector measure.
	\begin{defn}
		The twisted convolution of $\mu \in M(G \times \widehat{G})$ and a bounded $\nu \in M(G \times \widehat{G};X)$ is defined by 
		$$\mu \times \nu (A)=\int_{G  \times \widehat{G}} \int_{G  \times \widehat{G}} \chi_A\left(xx',\chi\chi'\right)  \chi'(x) d \mu(x,\chi)d \nu(x',\chi '),A \in \mathfrak{B} ( G \times \widehat{G}).$$
		whenever the mapping $(x',\chi')\mapsto\int_{G  \times \widehat{G}} \chi_A\left(xx',\chi\chi'\right)\chi'(x)\ d\mu(x,\chi)$ belongs to $L^1(G\times\widehat{G}, \nu).$
	\end{defn}
	Since $|\chi'(x)|=1,$ we have $\mu \times \nu \in M(G \times \widehat{G};X)$ and $\|\mu \times \nu \| \leq \|\mu\|\|\nu\|$.
	\begin{rem}
		It is worth noticing here that if the Banach space $X$ is $\mathbb{C},$ and $f,g \in L^1(G \times \widehat{G}) $ are considered as elements of $M(G \times \widehat{G})$, then $\mu_f \times \mu = \mu_{f \times g}$ where $d \mu_f(x,\chi)=f(x,\chi)d m_{G\times\widehat{G}}(x, \chi)$.
	\end{rem}
	For $f \in L^1(G \times \widehat{G})$ and $\nu \in M(G \times \widehat{G};X)$, define $f \times \nu=\mu_f \times \nu $ where $d \mu_f=fdm_{G\times\widehat{G}}.$
	We say that $f \times \nu \in C_c(G \times \widehat{G};X)$ if $d(f \times \nu)= f_\nu d m_{G\times\widehat{G}}, $ for some $f_\nu  \in C_c(G \times \widehat{G};X)$. The next proposition asserts the existence of such a function.
	\begin{prop} \label{P1}
		Let $\nu \in M(G \times \widehat{G}; X)$ be bounded. If $f \in C_c(G \times \widehat{G})$, then $f \times \nu \in C_c(G \times \widehat{G};X) $ and 
		$$\|f \times \nu\|_{C_c(G \times \widehat{G};X)} \leq \|f\|_{C_c(G \times \widehat{G})}\|\nu\|.$$
	\end{prop}
	\begin{proof}
		For $(x,\chi) \in G  \times \widehat{G}$ , we define $\textbf{f}_\nu(x,\chi) \in X $ by $$\textbf{f}_\nu(x,\chi):= \int_{G \times \widehat{G}} \tilde{f}_{(x,\chi)} d \nu$$ where $\tilde{f}_{(x,\chi)}(x',\chi')= f\left(x x'^{-1}, \chi \overline{\chi'}\right) \chi'(xx'^{-1})$. Now,
		\begin{align*}
			& \left\langle \mu_f \times \nu (A),x^* \right\rangle \\ =&  \int_{G  \times \widehat{G}} \int_{G  \times \widehat{G}} \chi_A\left(xx',\chi\chi'\right) \chi'(x) d \mu_f(x,\chi)d \langle \nu , x^* \rangle (x',\chi ') \\
			=&\int_{G  \times \widehat{G}} \int_{G  \times \widehat{G}} \chi_A\left(xx',\chi\chi'\right) \chi'(x) f(x,\chi)d m_{G\times\widehat{G}}(x,\chi)d  \langle \nu , x^* \rangle(x',\chi ')\\
			=& \int_{G  \times \widehat{G}} \int_{G  \times \widehat{G}} \chi_A\left(x,\chi\right) \chi'(x)\overline{\chi'(x')} f(xx'^{-1},\chi\overline{\chi})d m_{G\times\widehat{G}}(x,\chi)d  \langle \nu , x^* \rangle(x',\chi ')\\
			=&  \int_{G  \times \widehat{G}} \int_{A}f(xx'^{-1},\chi\overline{\chi}) \chi'(x) \overline{\chi'(x')} d m_{G\times\widehat{G}}(x,\chi)d  \langle \nu , x^* \rangle(x',\chi ') \\
			=& \int_A \int_{G \times \widehat{G}} \tilde{f}_{(x,\chi)}(x',\chi') d \langle \nu, x^* \rangle (x',\chi')d m_{G\times\widehat{G}}(x,\chi) \\
			=&  \left \langle \int_A \textbf{f}_\nu dm_{G\times\widehat{G}}  , x^* \right \rangle.
		\end{align*}
		Also,
		\begin{align*}
			\|\textbf{f}_\nu(x,\chi)-\textbf{f}_\nu(x',\chi')\|=&\left\|I_\nu\left(\tilde{f}_{(x,\chi)}-\tilde{f}_{(x',\chi')}\right)\right\| \leq \|\nu\| \left\|\tilde{f}_{(x,\chi)}-\tilde{f}_{(x',\chi')}\right\|_{C_c(G \times \widehat{G})}.
		\end{align*}
		Since $f \in C_c(G \times\widehat{G})$, therefore $ \textbf{f}_\nu \in C(G \times \widehat{G};X).$ Hence $d (f \times \nu) = \textbf{f}_\nu d m_{G\times\widehat{G}} $. Moreover,
		\begin{align*}
			\sup_{(x,\chi)\in G \times \widehat{G}} \|\textbf{f}_\nu(x,\chi)\|=&\sup_{(x,\chi) \in G \times \widehat{G}}\sup_{x^* \in B_{X^*}}\left| \left \langle I_\nu \left(\tilde{f}_{(x,\chi)}\right), x^*\right \rangle \right| \leq \|f\|_{C_c(G \times\widehat{G})}\|\nu\|. \qedhere
		\end{align*}	 
	\end{proof}
	Here is an analogue of Proposition \ref{TCProp1}.
	\begin{thm}
		Let $\nu \in M(G \times \widehat{G}; X)$ be bounded. Then for $1\leq p<\infty,$ if $f\in L^p(G \times \widehat{G})$, then $f\times\nu\in P_p(G \times \widehat{G};X)$ and $\|f \times \nu \|_{P_p(G \times \widehat{G};X)}\leq \|f\|_p \|\nu\|.$
	\end{thm}
	\begin{proof}
		From Proposition \ref{P1}, if $f \in C_c(G \times  \widehat{G})$, then $d( f \times \nu )=\textbf{f}_\nu dm_{G\times\widehat{G}}$ with $\textbf{f}_\nu \in C_c(G\times  \widehat{G};X)$. Then, for $x^* \in B_{X^*}$,
		\begin{align*} & \underset{G\times  \widehat{G}}{\int} |\langle \textbf{f}_\nu(x,\chi),x^* \rangle|^p d m_{G\times\widehat{G}}(x,\chi) \\
			 \leq & \underset{G\times  \widehat{G}}{\int} \bigg( \underset{G\times  \widehat{G}}{\int} |f(xx'^{-1},\chi\overline{\chi}) |d| \langle \nu, x^* \rangle | (x',\chi') \bigg) ^p dm_{G\times\widehat{G}}(x,\chi) \\
			\leq & \underset{G\times  \widehat{G}}{\int} \left(|\langle \nu, x^* \rangle|(G \times \widehat{G})\right)^{p-1}  \underset{G\times  \widehat{G}}{\int} |f(xx'^{-1},\chi \chi '^{-1})|^p d| \langle \nu, x^* \rangle | (x',\chi')   d m_{G\times\widehat{G}} (x,\chi) \\
			=& \left(| \langle \nu, x^* \rangle | (G \times  \widehat{G})\right)^{p-1} \underset{G\times  \widehat{G}}{\int}  \underset{G\times  \widehat{G}}{\int} |f(xx'^{-1}, \chi \overline{\chi'} )|^p dm_{G\times\widehat{G}}(x,\chi)d| \langle \nu, x^* \rangle | (x',\chi')  \\  
			=& \|f\|^p_{p} (|\langle \nu, x^*\rangle |(G\times  \widehat{G}))^p \leq  \|f\|^p_{p} \|\nu\|.
		\end{align*}
		Thus, the theorem holds true for functions in $C_c(G \times \widehat{G})$. Now, for $f \in L^p(G \times  \widehat{G})$, let $(f_n)$ in $C_c(G\times  \widehat{G})$ such that $\|f-f_n\|_p \to 0.$ Let $\nu_n=f_n \times \nu$. Then  $$ \left\|\nu_{n}-\nu_{m}\right\|_{p, m_{G\times\widehat{G}}}=\left\|\mathbf{f}_{\nu_{n}}-\mathbf{f}_{\nu_{m}}\right\|_{P_{p}(G \times \widehat{G};X)} \leq\left\|f_{n}-f_{m}\right\|_{L^{p}(G \times \widehat{G})}\|\nu\|.
		$$
		This shows that $\nu_{n}$ is a Cauchy sequence in $\mathcal{M}_{p}(G \times \widehat{G};X)$ which converges to $f \times \nu$ in $\mathcal{M}(G\times \widehat{G};X)$. Therefore, $f \times \nu \in P_{p}(\gtg;X)$ and
		\begin{align*}
			\|f \times \nu\|_{P_{p}(G\times \widehat{G};X)} =& \lim _{n}\left\|f_{n} \times \nu\right\|_{P_{p}(G\times \widehat{G};X)} \\ \leq& \lim _{n}\left\|f_{n}\right\|_{L^{p}(G\times \widehat{G})}\|\nu\|=\|f\|_{L^{p}(G\times \widehat{G})}\|\nu\| . \qedhere
		\end{align*}
	\end{proof}
	We finally prove the main result of this section.  This is the vector analogue of Young's inequality.
	\begin{thm}\label{VVYI}
		Let $ \nu \in M_{ac}(G \times \widehat{G};X)$. Then for $1<p< \infty$, if $\nu \in  M_{p}(G\times \widehat{G};X)$ and $f \in L^{q}(G\times \widehat{G})$ with $\frac{1}{p}+\frac{1}{q}>1,$ then $f \times \nu \in P_{r}(\gtg;X)$
		where  $\frac{1}{p}+\frac{1}{q}=1+\frac{1}{r}.$ Further
		$$\|f \times \nu\|_{P_{r}(G \times \widehat{G};X)} \leq\|f\|_{q}\|\nu\|_{p, m_{G\times\widehat{G}}}.$$
	\end{thm}
	\begin{proof}
		Since $\nu \in M_{ac}(G \times \widehat{G};X),$ by \cite[Pg. 259, Theorem 1]{D}, there exists $h_{x^{*}} \in L^{p}(G \times \widehat{G})$ such that $d\left\langle \nu, x^{*}\right\rangle=h_{x^{*}} d m_{G\times\widehat{G}}.$	Let $f \in C_c(G \times \widehat{G})$ and $\left\|x^{*}\right\|=1$. Then
		\begin{align*}
			\left| \left\langle f \times \nu(x,\chi), x^{*}\right\rangle\right |=&\left|\int_{G \times \widehat{G}} f(xx'^{-1},\chi\overline{\chi})  \chi'(x) \overline{\chi'(x')}d  \langle \nu , x^* \rangle(x',\chi ')\right|\\=&\left|\int_{G \times \widehat{G}} f(xx'^{-1},\chi\overline{\chi})  \chi'(x) \overline{\chi'(x')} h_{x^*}(x',\chi')d m_{G\times\widehat{G}}(x',\chi ')\right|\\ \leq & \int_{G \times \widehat{G}} |T^t_{(x',\chi')}f(x,\chi)|| h_{x^*}(x',\chi')| dm_{G\times\widehat{G}}(x',\chi') \\=&|f| * |h_{x^*} |(x,\chi)\end{align*} where $f*g(x,\chi)=\int_{G \times \widehat{G}}f(xx'^{-1},\chi \overline{\chi'})g(x',\chi')dm_{G \times\widehat{G}}(x',\chi')$ is the usual
			convolution.
		Now, by using the classical Young's inequality, we have$$
		\begin{aligned}
		\left(\int_{G \times \widehat{G}}\left|\left\langle f \times \nu(x,\chi), x^{*}\right\rangle\right|^{r} d m_{G\times\widehat{G}}(x,\chi)\right)^{1/r} & \leq\left\||f| *|h_{x^{*}}|\right\|_{L^{r}(G \times \widehat{G})} \\
		& \leq\|f\|_{L^{q}(G \times \widehat{G})}\left\|h_{x^{*}}\right\|_{L^{p}(G \times \widehat{G})} \\
		& \leq\|f\|_{L^{q}(G \times \widehat{G})}\|\nu\|_{p, m_{G\times\widehat{G}}}.
		\end{aligned}
		$$
		The result for an arbitrary $f \in L^q(G \times \widehat{G})$ follows via the denseness of $C_c(G \times \widehat{G})$ in $L^q(G \times \widehat{G}).$
	\end{proof}

 \section*{Declarations}
\noindent \textbf{Author's contributions:} 		 Both the authors have contributed equally.\\
 \textbf{Conflict of Interest:} 		 The authors declare that they have no conflict of interest.\\
 \textbf{Funding:} The first named author wishes to thank Graduate Aptitude Test in Engineering, India and Indian Institute of Technology Delhi for its research fellowship.\\
\noindent\textbf{Ethical Conduct:} Not applicable.\\
\textbf{Data Availability Statements:} 	The authors confirm that the data supporting the findings of this study are available within the article.

	\bibliographystyle{acm}
\bibliography{final.bib} 
\end{document}